\documentclass[a4paper,12pt]{article}
\usepackage[latin1]{inputenc}
\usepackage{amsmath,amssymb,amsthm,cite,verbatim}
\usepackage[all]{xy}
\usepackage[bookmarks=true]{hyperref}
\usepackage{color}
\usepackage{graphicx}

 \pdfoutput=1
\title{Automorphism group of split  Cartan modular  curves}
\author{ \ Josep Gonz\'alez  \footnote{The  author is partially supported by DGI grant  MTM2012-34611.\newline \emph{Keywords}: \,Automorphisms of modular curves, non-split Cartan subgroup.
}}

\newtheorem{prop}{Proposition}
\newtheorem{lema}{Lemma}
\newtheorem{teo}{Theorem}
\newtheorem{rem}{Remark}

\theoremstyle{definition}

\theoremstyle{remark}

\numberwithin{equation}{section}

\newcommand{\Q}{\mathbb{Q}}

\newcommand{\Z}{\mathbb{Z}}

\newcommand{\F}{\mathbb{F}}

\newcommand{\C}{\mathbb{C}}

\newcommand{\Gal}{\mathrm{Gal}}

\newcommand{\SL}{\operatorname{SL}}
\newcommand{\End}{\operatorname{End}}
\newcommand{\Aut}{\operatorname{Aut}}

\newcommand{\Jac}{\operatorname{Jac}}

\newcommand{\new}{\operatorname{new}}

\newcommand{\New}{\operatorname{New}}

\newcommand{\gl}{\mathfrak l}

\newcommand{\cO}{{\mathcal O}}

\newcommand{\cS}{{\mathcal S}}

\date{}
\begin{document}
\maketitle

\begin{abstract}
\noindent
We determine the automorphism group of the split Cartan modular curves $X_{\operatorname{split}}(p)$ for all primes $p$.
\end{abstract}

\section{Introduction}
For a prime integer $p$, let $X_{\operatorname{split}}(p)$ be the modular curve defined over $\Q$ attached to the congruence subgroup of level $p$
$$
\Gamma_{\operatorname{split}}(p)=\left\{ \left (\begin{array}{cc}a &b\\c& d\end{array}\right)\in\SL_2(\Z)\colon b\equiv c\equiv 0 \text{ or }a\equiv d\equiv 0\right\}\,.
$$
It is well-known that  $X_{\operatorname{split}}(p)$ is isomorphic over $\Q$ to the modular curve $X_0^+(p^2):=X_0(p^2)/w_{p^2}$, where $w_{p^2}$ stands for the Fricke involution. The genus of this curve is positive  when $p\geq 11$ and, in this case, it is at least $2$.

The automorphism group of the modular curve $X_0(N)$ was determined, except for $N=63$, by Kenku and Momose in \cite{KM} and was completed by Elkies in \cite{elkies}. Later,  the automorphism group of the modular curve $X_0^+(p)=X_0(p)/w_p$ was determined by Baker and Hasegawa in \cite{BH}. In this article, we focus our attention on  the automorphism group of the  split Cartan modular curves $X_{\operatorname{split}}(p)$.  Our main result is the following.

\begin{teo}\label{teo}
Assume that the genus of $X_{\operatorname{split}}(p)$ is positive. Then,
$$\Aut(X_{\operatorname{split}}(p))=\Aut_{\Q}(X_{\operatorname{split}}(p))\simeq\left\{\begin{array}{cc} \{1\} & \text{if $p>11$,}\\[5 pt]
(\Z/2\Z)^2 & \text{if $p=11$.}\end{array}\right.
$$
\end{teo}
\section{General facts}

   We recall that, for a normalized  newform $f=\sum_{n\geq 1} a_n q^n\in S_2(\Gamma_1(N))^{\operatorname{new}}$ and  a Dirichlet character  $\chi$ of conductor  $N'$, the function $$f_\chi:=\sum_{n\geq 1} \chi(n)a_n q^n$$
   is a cusp form in $S_2(\Gamma_1(\operatorname{lcm} (N, {N'}^2))$  (cf. Proposition 3.1 de \cite{Atkin}).  Here, as usual, $q=e^{2\pi \,i\,z}$.
   Let  $f\otimes \chi$ denote the unique normalized newform with $q$-expansion
$\sum_{n\geq 1} b_n q^n$ that satisfies
 $b_{\ell}=\chi(\ell) a_{\ell}$ for all primes  $\ell\nmid N\cdot N'$. If $f\otimes \chi$ is a newform of level $M$  and $f_\chi$ has level $M'$, then   $M|M'$.  Moreover, if $M=M'$, then $f_\chi=f\otimes\chi$.

\vskip0.2 cm

We  restrict ourselves to the cusp forms in $S_2(\Gamma_0(N))$. Let   $\New_N$  denote the set of normalized newforms in   $S_2(\Gamma_0(N))^{\operatorname{new}}$. For $f\in\New_N$,   $\varepsilon(f)$ denotes  the eigenvalue of  $f$ under the action of the Fricke  involution $w_N$ and set $\New_N^+=\{f\in \New_N\colon \varepsilon(f)=1\}$.
 For a cusp form $f=\sum_n b_n q^n \in S_2(\Gamma_0(N))$ such that   $\Q(\{b_n\})$ is a number field,  $S_2(f)$ denotes the $\C$-vector space  of cusp forms spanned  by  $f$ and its  Galois conjugates.  In the particular case that  $f\in\New_N$,     $A_f$ stands for  the abelian variety attached to $f$ by Shimura. It is well-known that $A_f$ is a quotient of $J_0(N):=\Jac (X_0(N))$  defined over  $\Q$ and the  pull-back of $\Omega^1_{A_f/\Q}$  is the $\Q$-vector subspace of elements in  $S_2(f) dq/q$ with rational $q$-expansion, i.e. $S_2(f) dq/q\cap\Q[[q]]$. Moreover, the endomorphism algebra $\End_\Q(A_f)\otimes \Q$ is a totally real number field.

\vskip0.2 cm
 From now on, we assume  $p\geq 11$ and   $\chi$  denotes the  quadratic Dirichlet character of conductor $p$, i.e. the Dirichlet character  attached to the quadratic number field  $K=\Q(\sqrt{p^*})$, where $p^*=(-1)^{(p-1)/2}$.
\begin{lema} The map $f\mapsto f\otimes \chi$ is a permutation of the set $\New_{p^2}\cup\New_p$. Under this bijection, there is a unique newform $f$, up to Galois conjugation, such that $f=f\otimes \chi$ when $p\equiv 3 \pmod 4$.

\end{lema}
\noindent{\bf Proof.}  Since    $f_\chi\in S_2(\Gamma_0(p^2))$ for $f\in\New_{p^2}\cup\New_p$ (cf. Theorem 3.1 of \cite{Atkin}),   the level of $f\otimes\chi$ divides $p^2$ and, thus, the map  is well defined. The injectivity follows from the fact that $(f\otimes\chi)\otimes\chi=f$.
The condition $f=f\otimes\chi$  amounts to saying that $f$ has complex multiplication (CM)  by the imaginary quadratic field $K$ attached to $\chi$, i.e. $p\equiv 3 \pmod 4$,  and, moreover, $f$ is obtained from a Hecke character $\psi$ whose  conductor is the ideal of $K$ of norm $p$, which implies  $f\in\New_{p^2}$. Since $f$ has trivial Nebentypus, the Hecke character $\psi$ is unique up to Galois conjugation.
\hfill $\Box$

\begin{rem} The above map does not preserve the eigenvalue of the corresponding Fricke involution, i.e.   it may be that  $\varepsilon(f)$ and $\varepsilon(f\otimes \chi)$ are different.
\end{rem}

\begin{rem} Let $f \in \New_{p^2}\cup \New_p$ without CM. If $f$ has an inner twist $\chi'\neq 1$, i.e. $f\otimes \chi'={}^{\sigma} f$ for some $\sigma\in G_\Q:=\Gal (\overline{\Q}/\Q)$, then $\chi'=\chi$ because $\chi'$ must be a quadratic character of conductor dividing $p^2$.
In such a  case,  $\End (A_f)\otimes\Q=\End_K(A_f)\otimes \Q$ is a non commutative algebra. Otherwise, $\End (A_f)\otimes\Q=\End_\Q(A_f)\otimes \Q$ is a totally real umber field.

\end{rem}
\begin{rem} If $f\in\New_{p^2}$ has CM, then the dimension of $A_f$ is the class   number of $K$, $A_f$ has all its endomorphisms defined over the Hilbert class field of $K$ and $\End_K(A_f)\otimes\Q$ is the CM field $\End_\Q(A_f)\otimes K$ which only contains the roots of the unity $\pm 1$ (cf. Theorem 1.2 of  of  \cite{gola011} and part (3) in Proposition 3.2 of \cite{Yang2}). Moreover, $f\in\New_{p^2}^+$ if, and only if, $p\equiv 3\pmod 8$ (cf. Corollary 6.3 of \cite{Yang1}).

\end{rem}

\begin{rem} For two distinct $f_1, f_2\in(\New_{p^2}\cup\New_p)/G_\Q$, the abelian varieties $A_{f_1}$ and $A_{f_2}$ are nonisogenous over $\Q$ and are isogenous if, and only if, $f_1\otimes \chi={}^{\sigma} f_2$ for some $\sigma\in G_\Q$ (see Proposition 4.2 of \cite{goji}) and, in this particular case, there is an isogeny defined over $K$.
\end{rem}

  The abelian variety $J_0^+(p^2):=\Jac (X_0^+(p^2))$ splits over $\Q$ as the product
$(J_0(p^2)^{\new})^{\langle w_{p^2}\rangle}\times J_0(p)$. More precisely,
\begin{equation}\label{split}
 J_0^+(p^2)\stackrel{\Q}\sim \prod_{f\in(\New_{p^2}^+\cup\New_p)/G_{\Q}} A_f\,.
\end{equation}

Each $f\in \New_p$ provides a vector subspace of $S_2(\Gamma_0(p^2))^{\operatorname{old}}$ of dimension $2$ genera\-ted by  $f(q)$ and $f(q^p)$. The normalized cusp forms $f(q)+p\,\varepsilon(f) f(q^p)$ and $f(q)-p\,\varepsilon(f) f(q^p)$ are eigenforms  for all Hecke operators $T_m$ with $p\nmid m$ and the Fricke involution $w_{p^2}$ with eigenvalues $1$ and $-1$ respectively.
 The splitting of $J_0^+(p^2)$  over $\Q$ provides the following decomposition for its vector space of regular differentials
\begin{equation}\label{diferencial}
\Omega^1_{J_0^+(p^2)}=\left(\bigoplus_{f\in\New_{p^2}^+/G_\Q}S_2(f(q))\frac{dq}{q}\right)\bigoplus\left(\bigoplus_{f\in\New_{p}/G_\Q} S_2(f(q)+p\,\varepsilon(f)f(q^p))\frac{dq}{q}\right)\,.
\end{equation}

 Let $g^+$ and $g_0$ be  the genus of the curves  $ X_0^+(p^2)$ and  $X_0(p)$,  respectively. From the genus formula for these curves, one obtains the following values
$$
\begin{array}{c|c|c|}
p  &   g^+ &g_0\\ \hline\hline
p\equiv 1 \pmod{12} &\displaystyle{\frac{(p-1)(p-7)}{24}}   & \displaystyle{\frac{p-13}{12}}\\[6 pt]
 p\equiv 5 \pmod{12} &\displaystyle{ \frac{(p-3)(p-5)}{24}  } & \displaystyle{\frac{p-5}{12} }\\[6 pt]
 p\equiv 7 \pmod{12}  & \displaystyle{ \frac{(p-1)(p-7)}{24}  }&\displaystyle{\frac{p-7}{12}} \\[6 pt]
 p\equiv 11 \pmod{12} &\displaystyle{ \frac{(p-3)(p-5)}{24} } & \displaystyle{\frac{p+1}{12}}\\\hline
\end{array}
$$

\section{Hyperelliptic case for $X_{\text{split}}(p)$}

\begin{prop}\label{11} Assume $p\geq 11 $. Then,
$X_{\operatorname{split}}(p)$ is hyperelliptic if, and only if,  $p=11$. Moreover, one has $$\Aut (X_{\operatorname{split}}(11))=\Aut_{\Q}(X_{\operatorname{split}}(11))\simeq
(\Z/2\Z)^2\,.
$$
\end{prop}
\noindent
{\bf Proof.} Assume  $X_0^+(p^2)$ is hyperelliptic. By  applying Lemma 3.25 of \cite{BGGP}, we obtain  $g^+\leq 10$, which implies  $p\leq 19$. We have $g^+=2$ if, and only if,   $p=11$ and, thus, the curve $X_0^+(11^2)$ is hyperelliptic. For $p>11$,  on has $g^+>2$   and, moreover, $p>13$ because  $X_0^+(13^2)$ is a plane quartic (cf. \cite{Baran}).  Lemma 2.5 of \cite{BGGP} states that there is a basis $f_1,\cdots, f_{g^+}$ of $S_2(\Gamma_0(p^2))^{\langle w_{p^2}\rangle}$ with rational $q$-expansions satisfying
\begin{equation}\label{congq}
f_i(q)=\left\{\begin{array}{cr} q^i+ O(q^i)&\text{ if the cusp $\infty$ is not a Weierstrass point of $X_0^+ (p^2)$,}\\q^{2 i-1}+ O(q^{2i-1})& \text{otherwise.}\end{array} \right.
\end{equation}
Moreover,  for any  such a basis,  the functions on $X_0 ^+(p^2)$ defined by
$$
x=\frac{f_{g^+}}{f_{g^+-1}}\,,\quad  y=\frac{ q dx/dq}{f_{g^+-1}}\,,
$$
satisfy $y^2=P(x)$ for  a unique squarefree polynomial $P(X)\in\Q[X]$ which has  degree $2g^++1$ or $2g^++2$ depending on whether $\infty$ is a Weierstrass point or not. The first part of the statement follows from the fact that,  for $p=17$ and $p=19$,  the vector space $S_2(\Gamma_0(p^2))^{\langle w_{p^2}\rangle}$ does not have any bases as in (\ref{congq}).

Now, we consider  $p=11$. In this case, $|\New_{11^2}|=|\New_{11}|=1$. Let $f_1\in\New_{11^2}^+$  and let $f_2\in\New_{11}$. The newform $f_1$ is the one attached to the elliptic curve $E_1/\Q$ of conductor $11^2$ with CM by  $\Z[(1+\sqrt{-11})/2]$, and $f_2$ is the newform attached to an elliptic curve $E_2/\Q$ of conductor $11$ without CM.  Since $\varepsilon(f_2)=-1$, the cusp forms  $f_1(q)$ and $h(q)=f_2(q)-11 f_2(q^{11})$  are a basis of $S_2(\Gamma_0(11^2))^{\langle w_{11^2}\rangle}$. Take the following functions on $X_0^+(11^2)$
$$
x=\frac{h}{f_1}=1 - 2 q + 2 q^3 - 2 q^5 + O(q^5)\,,\quad y=-2 \frac{q\,dx/q}{f_1}=4 - 8 q^2 + 8 q^3 + 24 q^4 - 32 q^5+ O(q^5)\,.
$$
Using $q$-expansions, we get the following equation for $X_0^+(11^2)$:
\begin{equation}\label{eq}y^2= x^6 -7 x^4+11 x^2+11\,.
\end{equation}The maps $(x,y)\mapsto (\pm x,\pm y) $ provide a subgroup of $\Aut_\Q(X_0^+(11^2))$ isomorphic to $(\Z/2\Z)^2$. Since $\Aut_{\overline\Q}(X_0^+(11^2))$ is a finite subgroup of $\End(E_1)\times \End(E_2)\simeq \Z[(1+\sqrt{-11})/2]\times\Z$, we have that it must be  a subgroup of $(\Z/2\Z)^2$, which proves the second part of the statement.
\hfill $\Box$

\section{Preliminary lemmas}
For an abelian variety  $A$ defined over a number field, we say that a number field $L$ is the splitting field of $A$ if it is the smallest number field where $A$ has all its endomorphisms defined.   We recall that $K$ is the quadratic field $\Q( \sqrt{p^*})$, where $p^*=(-1)^{(p-1)/2} p$, and $\chi$ is the quadratic Dirichlet character attached to $K$ .

\begin{lema}\label{split} Let $L$ be the splitting field of $J_0^+(p^2)$. If $p\equiv 3\pmod 8$, then
$L$  is de Hilbert class field de $K$. For $p\not\equiv 3 \pmod 8$, $L=K$ if  there exists
 $f\in\New_{p^2}^+\cup \New_p$ such that $f\otimes\chi\in\New_{p^2}^+\cup \New_p$, otherwise $L=\Q$.
\end{lema}
\noindent {\bf Proof.}  On the one hand, for two distinct $f_1,f_2\in(\New_{p^2}^+\cup \New_p)/G_\Q$  without CM, $A_{f_1}$ and $A_{f_2}$ are isogenous if, and only if, $f_2$ is the Galois conjugate of $f_1\otimes \chi$  and, in this case, the isogeny is defined over $K$.

On the other hand, if  $f\in\New_{p^2}^+\cup\New_p$ does not have CM, then $f$ has at most $\chi$ as an inner twist. In this case, the splitting field of $A_f$ is $K$ or $\Q$ depending on whether $\chi$ is an inner twist of $f$ or not. If $f$ has CM, then the splitting field of $A_f$ is the Hilbert class field of $K$  and $A_f$ is unique.
  \hfill $\Box$
\begin{lema}\label{defi}  All automorphisms  of $ X_0^+(p^2)$ are defined over $K$.

\end{lema}
\noindent {\bf Proof.} By Lemma \ref{split}, we only have to consider the case $p\equiv 3 \pmod 8$ and, by Proposition \ref{11}, we can assume $p>11$. Let $g_c$ be the dimension of the  abelian variety $A_f$ with $f\in\New_{p^2}$ having complex multiplication. We know that $g_c$ is the class number of $K$. Due to the fact that $g_c=(2V-(p-1)/2)/3$, where $V$ is the number of quadratic residues modulo $p$ in the interval $[1,(p-1)/2]$, we have $g_c\leq (p-1)/6$.  Since $g^+> 1+(p-1)/3$ for $p>11$ and $p\equiv 3 \pmod 8$, $g^+>  1+2 g_c$. Now, the statement is obtained  by applying  the same argument   used in the proof of Lemma 1.4 of \cite{KM}.
\hfill $\Box$

\begin{lema}
 The group $ \Aut_\Q (X_0^+(p^2))$ is isomorphic to $(\Z/2\Z)^s$ for some integer $s\geq 0$.
 \end{lema}
 \noindent {\bf Proof.}
Since  $\End_\Q(J_0^+(p^2))\otimes \Q$ is a product of totally real fields, any $u\in \Aut_\Q(X_0^+(p^2))$ acts as the identity or the product by $-1$ on $S_2(f)$ for $f\in\New_{p^2}^+$ and on $S_2(f(q)+p\varepsilon(f)f(q^p))$ for $f\in\New_p$.
 Hence, $\Aut_\Q(X_0^+(p^2))$ is isomorphic to a subgroup of $(\Z/2\Z)^r$ for $r$ equal to  $ |\New_{p^2}^+/G_\Q|+| \New_{p}/G_\Q|$. \hfill $\Box$

 \begin{lema}\label{overQ}
If  $u$ is  a nontrivial automorphism of $X_0^+(p^2)$, then  $u(\infty)$ is not a cusp.
\end{lema}
\noindent
{\bf Proof.}
 The modular curve $X_0(p^2)$ has $p+1$ cusps. Only the cusps $\infty$ and $0$ are defined over $\Q$. The remaining cusps $1/p,\cdots, (p-1)/p$ are defined over the $p$-th cyclotomic field $\Q(\zeta_p)$. The Galois group $\Gal(\Q(\zeta_p)/\Q)$ acts transitively on this set, and the cusp $w_{p^2}(i/p)=(p-i)/p$ is the  complex conjugate of the cusp $i/p$ (cf. \cite{Ogg}). Hence,  among the $(p+1)/2$  cusps of  $X_0^+(p^2)$  the cusp $\infty$ is the only one  defined over the quadratic field  $K$.

Let $u$ be  an  automorphism of $X_0^+(p^2)$ such that $u(\infty)$ is a cusp.   By Lemma \ref{defi},   $u(\infty)$ is defined over $K$ and, thus, $u(\infty)=\infty$.
Let $\Q [[q]]$ be the completion of the local ring $\cO_{X_0^+(p^2),\infty}$ with respect to the local parameter $q$. If $u(\infty)$ is the cusp  $\infty$, then
$u^*(q)=\zeta_m q+\sum_{n>1} c_i q^i$, where $\zeta_m$ is a primitive $m$-th root of unity and, thus,  $u$ is defined over $\Q (\zeta_m)$. Since the unique roots of  unity in $K$ are $\pm 1$, it follows that $u$ is defined over $\Q$.

Now,  we will prove that, if $u\in\Aut_\Q(X_0^+(p^2))$ is nontrivial, then $u(\infty)\neq \infty$.  For the hyperelliptic case $p=11$, the cusp $\infty$ has $(1,4)$ as $(x,y)$ coordinates in the equation given in (\ref{eq}). Hence, $\infty$ is not a fixed point for any of the three nontrivial involutions  of $X_{0}^+(11^2)$. Let  $p>11$. Since  $X_0^+(p^2)$ is nonhyperelliptic and, by Lemma \ref{overQ}, $u$ is an involution, we can exclude the case where all eigenvalues of $u$ acting on $\Omega^1_{J_0^+(p^2)}$ are equal to $-1$ and, thus,   $u$ must have eigenvalues $1$ and $-1$ acting on this vector space.
The vector space of cusp forms $\Omega^1_{J_0^+(p^2)} {q}/{dq}$ has a basis of normalized eigenforms (see (\ref{diferencial})). Since $u$ is defined over $\Q$, $u$ commutes with the Hecke operators and, thus, there are two normalized eigenforms such that their corresponding regular differentials $\omega_1=(1+\sum_{n>1} a_n q^n)dq$  and $\omega_2=(1+\sum_{n>1} b_n q^n)dq$ satisfy $u^*(\omega_1)=\omega_1$ and $u^*(\omega_2)=-\omega_2$. Hence, $u$ sends $\omega_1+\omega_2$, which does vanish at $\infty$, to $\omega_1-\omega_2$, which vanishes at $\infty$.
\hfill $\Box$

\vskip 0.2 cm

Let $\nu\in\Gal (K/\Q)$ be the conjugation corresponding to the Frobenius element of the prime $2$.
The following lemma is an adapted version of Lemma 3.3 of \cite{BH} to our context, restricted to the prime $2$.

\begin{lema}\label{gonality}
Assume that  $u\in\Aut (X_0^+(p^2))$ is nontrivial.  For any noncuspidal point $S\in X_0^+(p^2)({\C})$, the divisor
$$D_S:=(u T_2 -T_2 u^\nu)\left(\infty-S\right)\,,$$
where $T_2$ denotes de Hecke operator  viewed as a  correspondence of the curve $X_0^+(p^2)$, is nonzero but linearly
equivalent to zero.
In particular,  the $K$-gonality of $X_0^+(p^2)$ is at most $6$ and $u$ has at most $12$ fixed points.
\end{lema}
\noindent
{\bf Proof.}
Following the arguments used in Lemma 2.6 of \cite{KM}, one has $u\, T_2=T_2\, u^\nu$  and, thus, $D_S$ is a  principal divisor. Set $Q=u(\infty)$ and let  $P\in X_0(p^2)(\overline\Q)$ be such that $\pi^+(P)=Q$ where $\pi^+\colon X_0(p^2)\rightarrow X_0^+(p^2)$ is the natural projection. Since $Q$ is not a cusp, there is
an elliptic curve $E$  defined over $\overline\Q$ and a $p^2$-cyclic subgroup $C$ of $E(\overline\Q)$ such that $P=(E,C)$. The other preimage of $Q$ under $\pi^+$ is the point $w_{p^2}(P)=(E/C, E[p^2]/C)$.

  If $D_S$ is a  zero divisor, then  $u T_2 (\infty)$ must be equal to $T_2 u^\nu(\infty)$ because $T_2(\infty)=3\infty$ and $\infty$ is not in the support of $T_2(S)$.  To prove that $D_S$ is a nonzero divisor,  we only need to prove that the condition $3(Q)=T_2(Q^\nu)$ cannot occur.

Let $C_i$, $1\leq i\leq 3$, be the three $2$-cyclic subgroups of $E^\nu [2]$. Since $$T_2(Q^\nu)=\sum_ {i=1}^3 \pi^+((E^\nu/C_i,  C))\,,$$ the condition $3(Q)= T_2(Q^\sigma)$  implies that each elliptic curve $E^\nu/C_i$ is isomorphic to $E$ or  $E/C$. So, at least two quotients $E^\nu/C_i$ are isomorphic and, thus, $E^\nu$ must be an elliptic curve with CM by the  order $\cO$, where $\cO$ is the ring of integers  $\Z[\sqrt{-1}]$, $\Z[(1+\sqrt{-3})/2]$, $\Z[(1+\sqrt{-7})/2]$ or $\Z[(1+\sqrt{-15})/2]$. For the first three cases, $E^\nu$ is defined over $\Q$ and, in the last case, it is defined over $\Q(\sqrt 5)$. In any case, $E^\nu=E$. We  prove that, in all these case, there is a subgroup $C_i$ such that $E/C_i$ has CM by the order  $2\cO$  and this order contains an element of norm $2p^2$, which is not possible.

In the first case,  all elliptic curves $E/C_i$ are isomorphic and have CM by the order $2\cO$. Hence, $E/C_i$ must be  isomorphic to $
E/C$  and, thus, there is a $2p^2$-cyclic isogeny from  $E/C_i$ to itself.
This leads to the existence of an element in the order $2\cO$ of norm $2p^2$.

 In the second and third cases, $E$ is isomorphic to two curves $E/C_1$ and $E/C_2$, while $E/C_3$  is the elliptic curve with CM by $2\cO$. Hence, $E/C_3\simeq E/C$ and, thus, we can derive the existence of an element in $2\cO$ of norm $2p^2$.

In the last case, $E$ is defined over $\Q(\sqrt 5)$. Two  quotients $E/C_1$ and $E/C_2$ are isomorphic to the nontrivial Galois conjugated curve of $E$, and the curve
$E/C_3$ has CM by the order $2\cO$. Hence, $E/C_3$ is not isomorphic to $E$ and should be isomorphic to $E/C$. Again, we obtain  the existence of an element in $2\cO$ of norm  $2 p^2$.

By taking $S=u(\infty)$, $D_S$ is defined over $K$ and, thus, the $K$-gonality is at most $6$. Finally,  since $u^*(D_S)\neq D_S$ for some noncuspidal point $S\in X_0^+(p^2)({\C})$,   any  nontrivial automorphism of   $X_0^+(p^2)$ has at most $12$ fixed points (cf. Lemma 3.5 of \cite{BH}).
\hfill $\Box$
\begin{lema}
If $X_0^+(p^2)$ has a nontrivial automorphism and $p>11$, then $p\in\{ 17,19,23,29,31\}$.

\end{lema}
\noindent{\bf Proof.}
By applying Lemma 3.25 of \cite{BGGP} for the prime $2$, we obtain that $$g^+ < X_0^+(p^2)(\F_4)+1 \,.$$ By Lemma \ref{gonality}, the $\F_4$-gonality of $X_0^+(p)\otimes \F_{4}$ is $\leq 6$ and, thus,   $X_0^+(p)(\F_4)\leq 30$. Hence,  $g^+\leq 30$, which implies  $p\leq 31$. The algebra $\End(J_0^+(13^2))\otimes\Q $  is a totally real number field and, thus, it only contains the roots of unity $\pm 1$. Since  $X_0^+(13^2)$ is  nonhyperelliptic,   $\Aut(X_0^+(13^2))$  is  trivial and  we can discard the case $p=13$.
\hfill $\Box$

\begin{lema}\label{invo}
Every  nontrivial  automorphism  of $X_0^+(p^2)$ has  even order.
\end{lema}

\noindent {\bf Proof.}   Assume that there is  a nontrivial  automorphism $u$ of $X_0^+(p^2)$ whose order $m$ is odd. Let $X_u$ be the quotient curve $X_0^+(p^2)/ u$  and  denote by $ g_u$ its  genus.
Next, we  find a  positive lower bound $t$ for $g_u$.

 The endomorphism algebra $\End_K(J_0^+(p^2))\otimes \Q$ is the product
of some noncommutative algebras and some number fields $E_f=\End_K (A_f)\otimes\Q$ attached to the  newforms $f$ lying in a certain subset $\cS$ of $ (\New_{p^2}^+\cup \New_p)/G_\Q$. The set $\cS$ is formed by newforms $f$ without CM such that $f\otimes \chi\notin \New_{p^2}^+\cup\New_p$ ($E_f=\End_\Q(A_f)$ is a totally real number field) and by a newform $f$ with CM by $K$  if $p\equiv 3\pmod 8$ ($E_f=\End_\Q(A_f)\otimes K$ is a CM field). For $f\in\cS$,  the unique root of unities contained in $E_f$  are $\pm 1$. Since $m$ is odd,  the automorphism $u$ must act on each $E_f$ as the identity and, thus,  we have
$$t:= \sum_{f\in \cS}\dim A_f\leq g_u\,.$$
 An easy  computation provides the following values for $t$:
$$
\begin{array}{c|ccccc|}
p& 17 &19 &23 & 29& 31\\\hline\hline
g^+&7 & 9& 15& 26 & 30\\
t&5 & 5& 7& 15 & 12\\ \hline
\end{array}
$$
 Applying Riemann-Hurwitz formula,
 $$\displaystyle{m\leq \frac{g^+-1}{g_u-1}\leq \frac{g^+-1}{t-1}}<3\,,$$which yields a contradiction.
\hfill $\Box$
\section{Proof of Theorem \ref{teo}}

Assume that, for $p\in\{17,19,23,29,31\}$, there is a nontrivial  automorphism $u\in\Aut_K(X_0^+(p^2))$. By Lemma \ref{invo}, we can suppose that $u$ is an involution. Let  $g_u$ be the genus of the quotient curve $X_0^+(p^2)/u$.
   We  know that $u$   has at most $12$ fixed points.
 By Riemann-Hurwitz formula,  we get that the number of fixed points by $u$ must be even, say $2r$, and, moreover,
$$
g_u=\frac{g^++1-r}{2}\,,\quad 0\leq r\leq 6\,.
$$
If $g^+$ is even, then $u$ can have $2$, $6$ or $10$ ramification points, while for the case $g^+$ odd, $u$ can have $0$,  $4$, $8$ or $12$  such points.

For a prime $\ell\neq p$, the curve $X=X_0^+(p^2)$ has good reduction at $\ell$. Let  $\widetilde{X}$ be the reduction of  $X$  modulo $\ell$. We write
$$
N_{\ell}(n):=1+\ell^n-\sum_{i=1}^{2 g^+}\alpha_i^n\,,
$$
where $\alpha_1,\cdots,\alpha_{2 g^+}$ are the roots of polynomial
$$
\prod_{f\in \New_{p^2}^+\cup\New_p}(x^2-a_{\ell}(f) x+\ell)
$$
and $a_{\ell}(f)$ is the $\ell$-th Fourier coefficient of $f$. By Eichler-Shimura congruence,
$N_{\ell}(n)=|\widetilde{X}(\F_{\ell^n})|$.

Let $\gl$ be a prime of $K$ over $\ell$ with residue degree $s$. The reduction of $X\otimes K$ modulo $\gl$ is $\widetilde X\otimes \F_{\ell^s}$ which has an involution, say $\widetilde u$, with at most $2r$ fixed points.  The automorphism $\widetilde u$ acts on the set $\widetilde{X}(\F_{\ell^{s\, n}})$ as a permutation. If  $Q\in \cup_{i=1}^n\widetilde X(\F_{\ell^{ s\,i}})$, then the set $\cS_Q=\{ \widetilde{u}^i(Q)\colon 1\leq i \leq 2\}$ is contained in $ \cup_{i=1}^n\widetilde  X(\F_{\ell^{ s\,i}})$ and its  cardinality  is equal to $1$ or $2$ according to $Q$ is a fixed point of  $\widetilde u$ or not. Hence,    almost all integers  $R_{\ell}(n):=|\cup_{i=1}^n\widetilde  X(\F_{\ell^{ s\,i}})|$, $n\geq 1$,  are equivalent to the number of fixed points of $\widetilde{u}$ mod $2$
 and, moreover, the sequence $\{R_\ell (n)\}_{n\geq 1}$    can only contain at  most   $2r$ or $2r-1$ changes of parity  depending on whether $N_{\ell}(s)$ is even or odd. In other words, the sequence of integers  $\{ P_\ell(n)\}_{n\geq 1}$ defined by
$$0\leq P_\ell (n)\leq 1\,\,\text{ and }\,\,  P_{\ell}(n)=R_{\ell}(n+1)-R_{\ell}(n)\pmod 2\,,$$
can only contain at most $2r$ or $2r-1$ ones according to  $N_{\ell^s}$ being even or odd.

Note that the integer  $R_{\ell}(n+1)-R_{\ell}(n)$ can be obtained from the sequence $\{ N_\ell (s\, n)\}$
by using
$$\widetilde X (\F_{\ell ^{s\,d_1}})\cap \widetilde X(\F_{\ell ^{s\, d_2}})=\widetilde X (\F_{\ell ^{s\,\gcd(d_1,d_2)}})\,,\quad\text{ and if }
d_1|d_2 \text{ then }\widetilde X (\F_{\ell ^{s\,d_1}})\cup \widetilde X (\F_{\ell ^{s\,d_2}})=\widetilde X(\F_{\ell^{s\,d_2}})\,.
$$
More precisely, if $\{p_1,\cdots,p_r\}$ is the set of primes dividing $n+1$ and we put $d_i=(n+1)/p_i$ for $1\leq i\leq r$, then
$$
R_{\ell}(n+1)-R_{\ell}(n)= P_\ell (s(n+1))-\sum_ {j=1}^r (-1)^{r-1} \sum_{1\leq i_1<\cdots<i_j\leq r} P_\ell (s \gcd(d_{i_1},\cdots,d_{i_j}))\,.$$

For the five possibilities for $p$, we have:

\vskip 0.2 cm
{ $p=31$:}  $g^+=30$,   $2r\leq 10$ and   $\ell=2$  splits in $K=\Q(\sqrt{-31})$. One has
$$
N_2(1)=9\,,\qquad  \sum_{n\leq 36}P_2(n)=10\,.
$$

\vskip 0.2 cm
{$p=29$:} $g^+=26$,  $2r\leq 10$ and  $\ell=2$ is inert in $K=\Q(\sqrt{29})$. One has
$$
N_2(2)=42\,,\qquad  \sum_{n\leq 42}P_2(n)=11\,.
$$

\vskip 0.2 cm
{ $p=23$:}  $g^+=15$,   $2r\leq 12$ and  $\ell=2$  splits in $K=\Q(\sqrt{-23})$. One has
$$
N_2(1)=8\,,\qquad \sum_{n\leq 38}P_2(n)=13\,.
$$

\vskip 0.2 cm
{ $p=19$:} $g^+=9$,  $2r\leq 12$ and   $\ell=2$  is inert  in $K=\Q(\sqrt{-19})$. One has
$$N_2(2)=22\,,\qquad \sum_{n\leq 46 } P_2(n)=13\,.
$$

\vskip 0.2 cm
{ $p=17$:}  $g^+=7$,  $2r\leq 12$ and  $\ell=2$  splits in $K=\Q(\sqrt{17})$. One has
$$
N_2(1)=6\,,\quad\quad\sum_{n\leq 46}P_2(n)=13\,.
$$
So, we can discard the five cases considered and the statement is proved.

\bibliography{X_ns(11)}{}
\bibliographystyle{alpha}

\vskip 0.4 cm

\begin{tabular}{l}
Josep Gonz\'alez\\
\texttt{josepg\,\footnotesize{$@$}\,ma4.upc.edu}\\
Departament de Matem\`atica Aplicada 4  \\
Universitat Polit\`ecnica de Catalunya  \\
EPSEVG, Avinguda V\'ictor Balaguer 1\\
08800 Vilanova i la Geltr\'u, Spain\\[20pt]
\end{tabular}

\end{document}